\documentclass[12pt,a4paper]{article}
\usepackage{amstext}
\usepackage{fancyhdr}
\usepackage{tikz}
\usepackage{amsfonts,graphicx,bezier, amssymb}
\usepackage{caption}
\usepackage{pgfcore}
\usepackage{amsmath,amsthm,amssymb,amsfonts,amscd, array}
\usepackage{default}
\newtheorem{defn}{Definition}[section]
 \newtheorem{prop}{Proposition}
  
  \newtheorem{thm}{Theorem}[section]
  \newtheorem{cor}[thm]{Corollary}
  \newcommand{\Z}{\mathbb Z}
  \newtheorem{exam}{Example}
\newenvironment{prf}{\noindent{\bf{Proof:}}~~}{\hfill\rule{1ex}{1ex}\vskip1.5ex}
\newcommand{\Q}{\mathbb Q}
\newcommand{\R}{\mathbb R}

\newcommand{\F}{\mathbb F}

\newcommand{\beqa}{\begin{eqnarray}}
\newcommand{\enqa}{\end{eqnarray}}
\newcommand{\beq}{\begin{eqnarray*}}
\newcommand{\enq}{\end{eqnarray*}}

\newtheorem{qn}{Question}[section]
\newtheorem{rem}{Remark}[section]

\begin{document}

\begin{center}
{\bf\Large Nilpotent elements control  the structure of a  module}
\end{center}

\vspace*{0.3cm}
\begin{center}

David Ssevviiri\\

\vspace*{0.3cm}
Department of Mathematics\\
Makerere University, P.O BOX 7062, Kampala Uganda\\
E-mail: ssevviiri@cns.mak.ac.ug, ssevviirid@yahoo.com
\end{center}

\begin{abstract} A relationship between nilpotency and primeness in a module is investigated.
Reduced modules are expressed as sums of prime modules. It is shown that presence  of nilpotent module elements inhibits a module from 
possessing good structural properties.
  A general form is given of an example used in literature  to
 distinguish: 1)  completely prime modules from prime modules, 2) classical prime modules from classical completely prime modules, 
 and 3)  a module which satisfies the complete radical formula  from one which is neither 2-primal nor satisfies the radical formula.
   \end{abstract}


\vspace*{0.4cm}

{\bf MSC 2010} Mathematics Subject Classification: 16D70, 16D60, 16S90

\section{Introduction}\label{sec1}

Primeness and nilpotency are closely related and well studied notions for rings. We give  instances that highlight this relationship.
In  a commutative ring, the set of all nilpotent elements 
coincides with the intersection of all its prime ideals - henceforth called the prime radical.  
 A popular class of rings, called  2-primal rings (first defined in \cite{Gary} and later studied in \cite{Lee, Greg1, Greg2, Greg3} among others),
 is defined by requiring that in a not necessarily commutative ring, the set of all nilpotent
 elements  coincides with the  prime radical. In an arbitrary ring, Levitzki showed that the set of all   strongly
 nilpotent elements coincides with the intersection of all   prime ideals, \cite[Theorem 2.6]{McConnell}. 
 The upper nil radical of a ring which is the sum of all 
its nil ideals is the intersection of all its strongly prime ($s$-prime for short)  ideals, see \cite[p. 173]{Rowen} and 
\cite[Proposition 2.6.7]{RowenV1}. 
Every completely prime ring  has no nonzero nilpotent elements and 
every prime ring has no nonzero nilpotent ideals. A ring is semisimple (and hence a direct sum of
  prime rings) if and only if it is left (right) artinian
 and semiprime. A ring $R$ is semiprime (i.e., $R$ has no nonzero nilpotent ideals) if and only if $R$  is 
 isomorphic to a subdirect sum of prime rings, \cite[Theorem 4.27]{McCoy}. A unital ring $R$ is reduced (i.e., has no nonzero nilpotent 
 elements) if and only if $R$ is a subdirect sum of domains (completely prime rings), \cite[Example 3.8.16]{Gardner}.
 The upper nil radical of a ring $R$  is zero if and only if $R$ is a subdirect sum of strongly prime rings. Every ring which is a left essential
  extension of a reduced ring is a subdirect sum of rings which are essential extensions of domains, \cite{Beidar}.
  A semiprime ring $R$ has index $\leq n$ if and only if $R$ is a subdirect sum of prime rings of index $\leq n$, \cite{Efraim}.\\

 There is some effort to understand a   relationship between primeness and nilpotency in modules.
 Modules that  satisfy the radical formula (see \cite{Azizi, Jenkins, Man, Nikseresht, Pusat, Sharif} among others)
 were studied for this purpose. Behboodi in \cite{Behboodi1} and \cite{Behboodi2}  defined the
 Baer lower nil radical for modules  and sought its
 equivalence  with  the prime radical and the classical prime radical respectively. In \cite{SsevviiriNico} and its corrigendum,
 Groenewald and   I,  showed that in a uniserial module over a commutative ring, the set of all strongly nilpotent
 elements of a module coincides with the   classical prime radical of that module. Part of the aim  of 
  this paper is to continue with  this study that   establishes a relationship between nilpotency and primeness in a module. We   obtain
  structural theorems which relate reduced modules with prime modules, see Theorems \ref{VIT}, \ref{S}, \ref{Th1} \& \ref{zero}.
   It has been shown that absence of nilpotent module elements allows a module to behavave ``nicely'' by admitting certain structural
   properties, see paragraph after Question \ref{q}. The third major object of the paper is that, a general example has been 
   formulated for which a particular case was used in literature to distinguish several kinds of phenomena as we elaborate later.\\

   In Section 2, we give conditions under which a module is 
  nil. A link between nil left ideals and nil submodules is established leading to a possibility of using modules to prove K\"{o}the 
  conjecture in the negative. K\"{o}the conjecture which has existed since 1930 states that, the sum of two nil left ideals is nil, 
  see \cite[Theorem 2.6.35]{RowenV1}, \cite[p. 174]{Lam2} and \cite{Agata} for details about this conjecture.  \\
   
   In Section 3, we give 
   equivalent formulations for a module to be torsion-free and then for it to be reduced. We define co-reduced modules and 
   show that every reduced module over a commutative   ring is co-reduced. It is proved that 
   if $M$ is an injective module over a commutative  noetherian ring $R$ such that its indecomposable submodules are prime $R$-modules, then 
   $M$ is a reduced module which is   isomorphic to a direct sum of prime modules. 
   We introduce    $\mathcal{R}(M)$  the largest submodule of $M$ which is reduced as a module 
   and show that for a $\Z$-module $\Q/\Z$, 
   $\mathcal{R}(\Q/\Z)\cong \oplus_p~\Z/p\Z$ for all prime integers $p$.    
   This provides a general framework through which an already known result can be deduced, i.e., 
  a $\Z$-module $\Z/(p_1^{k_1}\times\cdots \times  p_n^{k_n})\Z$ is reduced if and only if $k_1=k_2=\cdots =k_n=1$.
 This in part answers Question 5.1 that was posed in \cite{S}.\\

    Section 4 is devoted towards getting conditions under which a reduced module is semisimple and hence a direct sum of prime  modules.
    We show that  if $M$  is  a reduced artinian $R$-module such that either $R$ or $R/\text{ann}_R(M)$ is embeddable 
     in $M$, 
  then every  $R/\text{ann}_R(M)$-module is  semisimple.  As corollaries, we get:  
  1) a faithful artinian completely prime module is semisimple; 
  2) an artinian, reduced and faithful free  module is semisimple; 
 3) an artinian, reduced, faithful and finitely generated module over a commutative ring is semisimple;
  4) a faithful holonomic module is reduced if and only if 
  it is semisimple and rigid; 5) an artinian $R$-algebra $M$  is a reduced $R/\text{ann}_R(M)$-module if and only if $M$ is a semisimple
   rigid $R/\text{ann}_R(M)$-module; and 6) if $(R, \mathfrak{m})$ is a Gorenstein artinian local ring, such that the injective envelope
   $\mathcal{E}_R(R/\mathfrak{m})$  of    $R/\mathfrak{m}$ is a reduced $R$-module, then $\mathcal{E}_R(R/\mathfrak{m})$
   is a semisimple $R$-module. It is also proved that, if $R$ is a left artinian ring and the zero submodule of an $R$-module
   $M$    satisfies the complete radical formula, then whenever $M$ is reduced, it follows that it is also semisimple. As a consequence,
   a module over a commutative artinian ring is reduced if and only if it is semisimple. This  retrieves  
    a well known result: a commutative artinian ring is reduced if and only if it is semisimple.\\

  In Section 5, examples are given to delimit and delineate the theory.
  For instance, we give what we call the Golden Example. It  serves several purposes. It is an 
  example of a prime nil module. It is a general form of the example which was used in \cite{S2} to 
  distinguish a prime module from a completely prime module, used in \cite{CCP} to distinguish a 
  classical prime module from a classical completely prime module and 
  an example used in \cite{cp2p} of a module which satisfies the complete radical formula but it is neither 2-primal nor  satisfies the radical
   formula. It is an example which shows another advantage (in addition  to those already known and given in \cite{S2}) that
   completely prime modules have over prime modules. If $M$  is a completely prime $R$-module, then $R/\text{ann}_R(M)$ is embeddable in $M$ and 
   if $M$ is in addition artinian, then $M$ as an $R/\text{ann}_R(M)$-module is semisimple. It is the same example that we have used to show that
    it is impossible to write a 2-primal module defined in \cite{2p} as one for which the prime radical  coincides  with a submodule generated by
     the nilpotent elements. This answers in negative Question 8.2.1 posed in \cite{thesis}.\\

   In Section    6, the last section, we give more effects of absence of nilpotent elements  
  on the primeness of a module. In particular, it is shown that a prime module  without nilpotent elements
   is $s$-prime, $l$-prime and completely prime.  If $M$ is a faithful $R$-module without nilpotent elements, then   
   $R$ has no nonzero nil ideals,  has no locally nilpotent ideals and it  has no nonzero   nilpotent ideals.\\
 
 A ring $R$ is {\it reduced} if it has no nonzero nilpotent elements.
 So, a ring with a nonzero nilpotent element is not reduced. This was the motivation in \cite{SsevviiriNico} for defining a nilpotent 
  element of a module. Let $R$ be a ring. 
  An $R$-module $M$  is {\it reduced} (see \cite{Baser, Leezhou, Rege})
  if for all $a\in R$ and $m\in M$, $a^2m=0$  implies that $aRm=0$. Reduced modules were called completely semiprime modules in \cite{CCP}.
  As for rings, we say   that an $R$-module is not reduced if it has a nonzero nilpotent element.
  
  \begin{defn}\rm
   An element $m$ of an $R$-module $M$ is {\it nilpotent} if either $m=0$  or there exists a positive integer $k$ and an element 
   $a\in R$ such that    $a^km=0$ and $aRm\not=0$. In this case, $a$ is called the {\it nilpotentiser} of $m$ and the least such $k$ is called
    the {\it degree of nilpotency} of $m$ with respect to $a$.
  \end{defn}
  
  \begin{defn}\rm 
 An $R$-module $M$ is {\it  nilpotent} if for each  $0\not=m\in M$, there exists\footnote{$a(m)$ means that $a$
  is dependent on $m$.}  $a(m)\in R$ and  a fixed $k\in \Z^+$  independent  of $m$   such that $aRm\not=0$ and $a^km=0$. 
 \end{defn}

 Throughout this paper, unless stated otherwise, $R$ will denote a unital  associative ring 
 and  $M$  will be a unital left $R$-module. We write 1 for the unity of  a ring $R$. By $R$-Mod, 
     $M_n(R)$ and $\text{End}_R(M)$   we shall respectively mean the category of all $R$-modules, the ring of all $n\times n$ 
     matrices over $R$ and   the endomorphism ring of $M$ over $R$. 
 $\text{ann}_R(M)$ (resp. $\text{ann}_R(m)$)  will denote the ideal $\{r\in R~:~rM=0\}$ (resp. the left ideal  $\{r\in R:rm=0\}$) of $R$.
 For an $R$-module $M$ and $r\in R$, $(0:_M r)$ will denote the set $\{m\in M~:~rm=0\}$ of $M$.
 $\mathcal{N}(M)$  denotes the collection of all nilpotent elements of a module $M$.

\section{Nil modules}

\begin{defn}\rm 
 A module is {\it nil} if all its elements are nilpotent.
\end{defn}

  It is clear that a nilpotent module  is nil.
  If the ring  $\text{End}_R(M)$ or $R$ is  nil\footnote{Note that, a nil ring cannot be unital - this is an exception to the 
  general rings used in this paper.} and for all nonzero elements $m\in M$, $\text{ann}_R(m)\subsetneq R$; then $M$ is nil. It 
    is also easy to see that, if the ring  $\text{End}_R(M)$ or $R$ is  nilpotent  and for all nonzero elements $m\in M$, 
    $\text{ann}_R(m)\subsetneq R$; then $M$ is nilpotent. Fischer \cite[Theorem 1.5]{Fisher} gave a condition for the ring $\text{End}_R(M)$
    to be nilpotent.\\
   
   A nil module need not be nilpotent. Let $T_n$ be a ring of all $n\times n$ matrices over a division ring  where $n\geq 2$. 
   Each module $_{T_n}T_n$ is nilpotent. Consider the  direct sum $A:= \oplus_{n=2}^{\infty}T_n$
   which has an ascending chain $$ T_2\subset T_2 \oplus T_3\subset \cdots \subset \oplus_{n=2}^{k}T_n\subset \cdots $$ of
    submodules such that $A$ is the union  $A= \cup_{n=2}^{\infty}(\oplus_{n=2}^k T_n)$ of the members of the chain. $A$ is nil $A$-module 
    which is not nilpotent.\\

  A left ideal $I$ of a ring $R$ is {\it dense} if given any $0\not=r_1\in R$ and $r_2\in R$, there exists $r\in R$ such that 
 $rr_1\not=0$ and $rr_2\in I$. An  ideal $I$ of a ring $R$ is {\it essential} if for all nonzero ideals $J$ of $R$, $I\cap J\not=0$.
 Every dense  ideal is an essential ideal.
 
 \begin{prop}\label{l} If   $R$ is a ring with  a left ideal which is both  nil and dense, then the module $_RR$ is nil.
  \end{prop}

\begin{prf} Let $J$ be a nil and dense left ideal of $R$. Suppose $0\not=m\in R$. $1\in R$ and by definition
   of dense left ideals,  there exists $r\in R$ such that  $rm\not=0$ and $r1=r\in J$. So, $rRm\not=0$.
   Since $J$ is nil, $r$ is nilpotent and     hence $r^km=0$ for some $k\in \Z^+$.
   \end{prf}

The envelope of a submodule $N$  of an $R$-module $M$ is the set 
 $$E_M(N):=\{rm ~:~r^km \in N, r\in R, m\in M, k\in \Z^+\}.$$ This set  was used to define modules that satisfy the radical formula,
 see \cite{Azizi, Nikseresht, Pusat, Sharif} among others. In the context of modules that satisfy the radical formula, 
 $E_M(0)$ was considered to be the module analogue of the set $\mathcal{N}(R)$ of all nilpotent elements of a ring $R$.  Since $E_M(0)$ is not 
 a submodule, we write the submodule of $M$  generated by $E_M(0)$ as $\langle E_M(0) \rangle$.

 \begin{prop}\label{nil}
  For any $R$-module $M$:
  \begin{enumerate}
   \item $\langle E_M(0)\rangle \subseteq \langle \mathcal{N}(M)\rangle $;
   \item if $E_M(0)=M$, then $\mathcal{N}(M)=M$ and hence $M$  is nil;
   \item if $m\in \mathcal{N}(M)$, then there exists $r\in R$ such that $rm\in E_M(0)$.
  \end{enumerate}

 \end{prop}

 \begin{prf}
    If $0\not=m\in E_M(0)$, then $m=rn$ and $r^kn=0$ for some $r\in R$,  $k\in \Z^+$ and $n\in M$. This shows that $n\in \mathcal{N}(M)$.
   So, $m=rn\in \langle \mathcal{N}(M)\rangle$ and $E_M(0)\subseteq \langle \mathcal{N}(M)\rangle$. Hence, 
   $\langle E_M(0) \rangle \subseteq \langle \mathcal{N}(M)\rangle$. 2 follows immediately from 1 and 3 is trivial.
  \end{prf}

   Let $R$ be a commutative ring.   If every nonzero endomorphism $f_r$ of an $R$-module    $M$ which is given by $f_r(m)=rm$ for 
 some $r\in R$ and $m\in M$ is both nilpotent and surjective, then    $M=E_M(0)$ and $M$  is nil. For if $f_r$ is both nilpotent and surjective, then for all $n\in M$, there exists $m\in M$ such that $n=rm$.
  Since $f_r$ is nilpotent,   there exists a positive integer $k$ such that  $r^km=0$. It follows that $n\in E_M(0)$ and hence
  $M=E_M(0)$. Nil modules obtained this way are secondary modules. A  module $M$  is {\it secondary } \cite{Macdonald} if $M\not=0$ and 
    for each $r\in R$, the endomorphism $m \mapsto rm$ of $M$ is either surjective or nilpotent.\\

  On the other hand, it is impossible for a  nilpotent endomorphism $f_r(m)=rm$ of $M$ to be injective and hence it cannot be 
  an isomorphism.  Suppose that $f_r(m)=rm\not=0$   and $f_r$ is nilpotent, i.e., $r^km=0$ for some $k\in \Z^+$. Define $m_1=rm$. 
  Then $f_r(m_1)=rm_1=r^2m$. Take $m_2=rm_1$, we get    $f_r(m_2)=rm_2= r^2m_1 =r^3m$, continuing this way, we get
  $f_r(m_{k-1})= r^km=0$. $f$ injective implies that $m_{k-1}=0$ so that $f(m_{k-2})=0$. 
  Continuing with this process leads to $m_1=rm=0$ which is a contradiction.\\



We here repeat  Question 8.2.2 posed in \cite{thesis}. 

 \begin{qn} [\rm{\bf Module analogue of K\"{o}the conjecture}]\label{kothe} 
  Is the sum of nil submodules nil?
 \end{qn}

   The importance of Question \ref{kothe} lies in the fact that it can be used to solve K\"{o}the conjecture in the negative. For if
   $I_1$ and $I_2$ are nil left ideals of a ring $R$ and $M$  is an $R$-module, then the submodules $I_1m$ and $I_2m$  of $M$ are also
    nil. If the sum $I_1m+I_2m =(I_1+I_2)m$ is not a nil submodule, then it follows that the sum $I_1+I_2$ is not a nil left ideal of
    $R$ and hence K\"{o}the conjecture would be false.\\
    
    It is tempting for one to think that may be $\langle \mathcal{N}(M)\rangle= \langle E_M(0)\rangle$ for any module $M$.
    However, for a $\Z$-module      $M:=\Z/ 4\Z$,  $\langle E_M(0)\rangle=0$ and $\langle \mathcal{N}(M)\rangle= M$.

   \begin{defn}\rm
    A submodule $P$ of an  $R$-module $M$    is {\it prime} if for all ideals $A$ of $R$ and submodules $N$  of $M$,
    $AN \subseteq P$  implies that either $AM \subseteq P$  or $N \subseteq P$. A module is prime if its zero submodule is a prime submodule.
   \end{defn}

    \begin{prop}\label{B}
     Let $\{a_i\}_{i\in I}$ be a collection of all nilpotentisers of elements of an $R$-module $M$  with corresponding degrees of nilpotency
     $\{k_i\}_{i\in I}$.  Then  $$\mathcal{N}(M) \subseteq \sum_{i\in J\subseteq I} (0:_{M} ~{a_i}^{k_i}).$$ We get equality whenever $M$ 
     is nil. 
     \end{prop} 

    \begin{prf}
     If $m_i\in\mathcal{N}(M)$, then there exists $a_i\in R$ and $k_i\in \Z^+$ such that $a_i^{k_i}m_i=0$ and $a_iRm_i\not=0$. This implies that
     $m_i\in (0:_M~ a_i^{k_i})$ and hence $\mathcal{N}(M)\subseteq \sum_{i\in J\subseteq I}(0:_M~a_i^{k_i})$. To see why $i\in J\subseteq I$, 
     if $m\in M$  has several nilpotentisers, then we take only one leading to $J\subseteq I$. If $M$  is nil, $\mathcal{N}(M)=M$ and  it follows that 
     $M=\sum_{i\in J\subseteq I} (0:_{M} ~{a_i}^{k_i})$.
     \end{prf}


 \begin{exam}\rm
  Let $R:= \F_2[x]/\langle x^3\rangle$. The module $_RR$ has all its elements nilpotent except $x^2$. It is easy to see that 
  $\mathcal{N}(_RR)\subsetneq  \sum_{i\in I } (0:_R a_i^{k_i})=R$.
 \end{exam}
 
 \begin{exam}\rm 
  If $R:= M_n(D)$ the matrix ring of order $n$ over a division ring $D$, then the module $_RR$ is nil.  
  So, $R=\mathcal{N}(_RR) =  \sum_{i\in J\subseteq I} (0:_R a_i^{k_i})$.




 \end{exam}

 \begin{exam}\rm
 Let $M:=M_n(\Z/k\Z)$ where $k$  is an integer be the ring of all $n\times n$  matrices defined over the ring $\Z/k\Z$,
 $$N:=\left\{(m_{ij})\in M~:~ \sum_{j=1}^nm_{ij}=0~(\text{mod}~k) ~\forall i\in \{1, 2, \cdots, n\}\right\}~\text{and}~ 
 R:=M_n(\Z).$$ $N$  is a   nil $R$-module and $N=\mathcal{N}(N)=\sum_{i\in J\subseteq I}(0:_N a_i^{k_i})$.
 \end{exam}

 \section{Reduced modules}
 
 An $R$-module $M$   is   {\it completely prime } \cite[Definition 2.1]{DS}    if for all elements 
 $a\in R$ and $m\in M$,   $am=0$  implies that either $aM=0$ or $m=0$.  An $R$-module $M$   is 
 {\it rigid} \cite{B} if for all  $m\in M$,  $a\in R$ and a   positive integer $k$,  $a^km=0$ implies that $am=0$. For modules over 
 commutative rings, rigid modules are indistinguishable from  reduced modules. \\
 
For a given ring $R$, it is easy to see that
$$\left\{\text{torsion-free $R$-modules}\right\}\subsetneq \left\{\text{completely prime $R$-modules}\right\}
\subsetneq \left\{\text{reduced $R$-modules}\right\}$$
$$\subsetneq \left\{\text{ rigid $R$-modules}\right\}.$$

In Proposition \ref{lab}, we give equivalent statements in terms of reduced modules and completely 
prime modules for a module to be torsion-free.
 
 \begin{prop}\label{lab}
  For an $R$-module $M$, the following statements are equivalent:
  \begin{enumerate}
   \item $M$  is torsion-free;
      \item for all nonzero elements $m\in M$, $\text{ann}_R(m)=0$;
   \item $M$ is a reduced  module and $\text{ann}_R(Rm)=0$ for all nonzero elements $m\in M$;
   \item $M$  is a completely prime  module and $\text{ann}_R(M)=0$.
  \end{enumerate}

 \end{prop}

\begin{cor}
 A completely prime module is faithful if and only if it is torsion-free.
\end{cor}

\begin{prop}\label{red}
 If $R$ is a commutative ring and $M$ is an $R$-module, then the following statements are equivalent:
 \begin{enumerate}
  \item $M$  is a reduced $R$-module;
  \item for every nonzero $m\in M$, the $R$-module $R/\text{ann}_R(m)$ is reduced;
  \item for every nonzero $m\in M$, the cyclic $R$-module $Rm$ is reduced;
  \item every minimal submodule of $M$  is reduced;
  \item every nonzero endomorphism $f$ of $M$,  of the form $f_a(m)=am$, where $a\in R$ and $m\in M$ is not nilpotent;
  \item for every nonzero $m\in M$, $\text{ann}_R(m)$ is a semiprime ideal of  $R$;
  \item for every $r\in R$ and $k\in \Z^+$, $\text{Ker}~f_r = \text{Ker}~f_r^k$, where $f_r(m)=rm$ with $m\in M$;
  \item $(0:_Mr)=(0:_Mr^k)$ for all $r\in R$ and $k\in \Z^+$;
  \item $E_M(0)=0$;
  \item $\mathcal{N}(M)=0$.
  
 \end{enumerate}
\end{prop}

Proposition \ref{red} makes it possible to dualise reduced modules defined over commutative rings.

\begin{defn}\rm Let $R$ be a commutative  ring.  An $R$-module $M$  is {\it co-reduced} if  for any  endomorphism 
 $f_r(m)=rm$ of $M$ with $r\in R$,  $$\text{Co-Ker}~f_r = \text{Co-Ker}~f_r^k.$$  
\end{defn}

Suppose $M$  is a module over a commutative  ring and for every $r\in R$, $\text{Ker}~f_r$ (resp. $\text{Im}~f_r$) is a maximal (resp. minimal) submodule of $M$ 
 where $f_r:M\rightarrow M$ is the endomorphism $f_r(m)=rm$, then $M$  is reduced (resp. co-reduced). Note that 
  $\text{Ker}~f_r\subseteq \text{Ker}~f^2_r\subseteq \text{Ker}~f^3_r\subseteq \cdots $  
  (resp. $\text{Im}~f_r\supseteq \text{Im}~f^2_r\supseteq \text{Im}~f^3_r\supseteq \cdots $) is an ascending chain (resp. descending chain)
  and is constant when  $\text{Ker}~f_r$ (resp. $\text{Im}~f_r$) is a maximal (resp. minimal) submodule of $M$.

  \begin{prop} If a module which is defined over a commutative ring is reduced and has a finite number of elements, then it is coreduced.
   
  \end{prop}
\begin{prf}
 Let $f_r:M\rightarrow M$  be given by $f_r(m)=rm$. Then for a positive integer $k$, $f_{r^k}(m)= r^km$. By the first Isomorphism Theorem,
 $M/(0:_M~r)\cong rM$ and $M/(0:_M~r^k)\cong r^kM$. If $M$  is reduced, by Proposition \ref{red},
 $(0:_M~r)=(0:_M~r^k)$. This implies that  $rM\cong r^kM$.
 However, $r^kM\subseteq rM$ and $M$  has finite order. It must follow that $r^kM=rM$ and hence $M$  is co-reduced.
\end{prf}

\begin{thm}\label{VIT}If $M$ is an injective module defined over a commutative noetherian ring  $R$ such that its indecomposable
submodules are prime $R$-modules, then $M$  is a reduced module which is  isomorphic to a direct sum of prime modules.
\end{thm}

\begin{prf}
 Any injective module $M$ over a commutative noetherian ring $R$   is isomorphic to $\oplus_{i} \mathcal{E}_R(R/P_i)$
 where $P_i$ are prime  ideals of $R$  and $\mathcal{E}_R(R/P_i)$ is the injective envelope of the $R$-module $R/P_i$.
 If each indecomposable component $\mathcal{E}_R(R/P_i)$ of $M$ is a prime $R$-module, then it is also reduced.
  However, a direct sum of reduced modules is reduced. This shows that $M$  is a reduced module which is isomorphic to 
 a direct sum of prime modules.
\end{prf}

\begin{exam}\rm
$\Q/\Z$ is an injective module over $\Z$ (a commutative  noetherian ring). However it is not reduced.  $\frac{1}{4}$ is a nilpotent element
in the $\Z$-module $\Q/\Z$. $2^2 \times \frac{1}{4}\in \Z$ but $2 \times \frac{1}{4}\not\in \Z$. We observe that
$$\Q/\Z\cong \oplus_{p}~\Z(p^{\infty})\cong \oplus_{p}~ \mathcal{E}_{\Z}(\Z/p\Z)$$ since by \cite[Example 3.36]{Lam}
every injective envelope of a $p$-group  is a prufer $p$-group $\Z(p^{\infty})$.
The indecomposable components $\mathcal{E}_{\Z}(\Z/p\Z)$ are not prime $\Z$-modules.
 \end{exam}
 
 \begin{exam}\rm 
 $\Q$ is an injective module over $\Z$ (a commutative  noetherian ring). $\Q=\mathcal{E}_{\Z}(\Z)$ which is a prime $\Z$-module and 
 hence a reduced $\Z$-module.
 \end{exam}

 For an $R$-module $M$,  let $$\mathcal{R}(M):=\{m\in M~:~ (a^km=0) \Rightarrow (aRm=0)~\text{for}~ a\in R~  \&~ k\in \Z^+\}.$$ $\mathcal{R}(M)$ is
  the collection of all non-nilpotent elements of the $R$-module $M$ together with the zero element.  A module  $M$ is reduced  if and only if 
  $\mathcal{R}(M)=M$ and $M$ is nil if and only if $\mathcal{R}(M)=0$. If $M$  is a module over a commutative ring $R$, then the set 
  $\mathcal{R}(M)$ is a submodule of $M$. For if $m_1, m_2\in \mathcal{R}(M)$ and there exists $a\in R$ such that  
  $a^k(m_1+m_2)=0$ but $aR(m_1+m_2)\not=0$. Then either $aRm_1\not=0$ or $aRm_2\not=0$. If $aRm_1\not=0$ and $aRm_2=0$. Then 
  $a^tm_1\not=0$ for all $t\in \Z^+$ since $m_1\in \mathcal{R}(M)$  and $a^tm_2=0$. It follows that $a^k(m_1+m_2)\not=0$  which is 
  a contradiction. Similarly, if $aRm_1=0$ and $aRm_2\not=0$, then $a^k(m_1+m_2)\not=0$  leading to a contradiction. Now suppose that 
  $aRm_1\not=0$ and $aRm_2\not=0$.  It follows that $a^tm_1\not=0$ and $a^tm_2\not=0$ for all $t\in \Z^+$ 
  since $m_1, m_2\in \mathcal{R}(M)$. This implies that
  the submodules $a^kRm_1$ and $a^kRm_2$ of $M$ are both nonzero and hence their sum $a^kR(m_1+m_2)$ is also nonzero. However, 
    $a^k(m_1+m_2)=0$   implies that $a^kR(m_1+m_2)=0$ which is also a contradiction. This shows that if $m_1, m_2\in \mathcal{R}(M)$,
  then $m_1+m_2\in \mathcal{R}(M)$. If $m\in \mathcal{R}(M)$ and $r\in R$ such that  $a^k(rm)=0$ for some $k\in \Z^+$ 
  and $a\in R$ but $aR(rm)\not=0$, then $(ar)Rm\not=0$. By definition of $\mathcal{R}(M)$,  $(ar)^tm\not=0$ for all $t\in \Z^+$ so that 
  $a^k(rm)\not=0$ which is a contradiction. This shows that $r\in R$ and $m\in \mathcal{R}(M)$  implies that $rm\in \mathcal{R}(M)$.\\
  
  For modules over commutative  rings, $\mathcal{R}(M)$ is the largest submodule of $M$ which is reduced as a module.

 \begin{thm}\label{S} As $\Z$-modules, 
   $$\mathcal{R}(\Q/\Z)\cong  \oplus_{p} ~\Z/p\Z$$ for all prime integers $p$. Hence, $\mathcal{R}(\Q/\Z)$ is a reduced  $\Z$-module
   which is isomorphic to a direct sum of prime $\Z$-modules.
  \end{thm}
  
  \begin{prf}
   $\mathcal{R}(\Q/\Z)= \mathcal{R}(\oplus_{p}~ \Z(p^{\infty}))=  \oplus (\left( \frac{1}{p}\Z\right)/\Z) \cong \oplus_{p}~ \Z/p\Z$,
   where $\left( \frac{1}{p}\Z\right)/\Z$ is the cyclic subgroup of $\Z(p^{\infty})$  with $p$ elements; it contains those 
   elements of $\Z(p^{\infty})$ whose order divides $p$ and corresponds to the set of $p$-th roots of unity.
    Each $\Z$-module $\Z/p\Z$ is a simple $\Z$-module and hence a prime $\Z$-module.
  \end{prf}

\begin{cor}\label{c1}
 The $\Z$-module $\Z/(p_1\times \cdots \times p_k)\Z$ is  reduced whenever the prime integers $p_i$ are all distinct.
\end{cor}

\begin{prf}
 $\Z/(p_1\times \cdots \times p_k)\Z\cong \Z/p_1\Z \oplus \cdots \oplus \Z/p_k\Z$ which is a submodule of a reduced $\Z$-module
 $\mathcal{R}(\Q/\Z)\cong  \oplus_{p} ~\Z/p\Z$.
\end{prf}

Theorem \ref{S}  partly answers in affirmative Question 5.1 in \cite{S} which asks whether  it is possible to obtain
a general framework through which   results like Corollary \ref{c1} which were obtained in \cite{S}  can be retrieved. \\

Unlike prime rings which are closed under essential extension, prime modules need not be closed under essential extension.
The $\Z$-module $\Z/p\Z$ is a prime $\Z$-module. However, its essential extension $\mathcal{E}_{\Z}(\Z/p\Z)$ which is the prufer $p$-group
$\Z(p^{\infty})$ is not a prime $\Z$-module.\\

\begin{table}[t]
\begin{center}

 \begin{tabular}{|llll|}\hline
   & $M$  &  $R$  &  $\mathcal{R}(M)$ \cr\hline
   1. & $\F_2[x]/\langle x^2\rangle$ &  $\F_2[x]/\langle x^2\rangle$ & $\{0, x\}$ \cr
   2. & $\F_2[x]/\langle x^3\rangle$ & $\F_2[x]/\langle x^3\rangle$ & $\{0, x^2\}$ \cr
   3. & $\Q/\Z$                      & $\Z$                         &   $\oplus (\left( \frac{1}{p}\Z\right)/\Z)$  \cr
      &                               &                              & $\cong  \oplus_{p} ~\Z/p\Z$ \cr 
   4. &  $\Z/(p_1^{k_1}\times \cdots \times p_n^{k_n})\Z$  & $\Z$  &  $\{0\}\cup\{mp_1^{k_1-1}\cdots p_n^{k_n-1}\}_{m=1}^{p_1p_2\cdots p_n -1}$  \cr
      &                                                   &      &   $\cong\Z/(p_1\times \cdots \times p_n)$                 \cr
   5. &  $\Z_p$, $p$-adic integers & $\Z$ &  $\prod_{i=1}^np_i^{t_i}\Z_{p_i}$, $t_i\in \Z^+$ \cr
   6. & $\R[x_1, \cdots, x_n]$     & $\R[x_1, \cdots, x_n]$ &  $\R[x_1, \cdots, x_n]$ \cr
   7. & $\R[x_1, \cdots, x_n]$     & $\R $ &  $\R[x_1, \cdots, x_n]$ \cr\hline
 \end{tabular}
 
 \caption{Examples of   $\mathcal{R}(M)$}\label{ti}
 \end{center}
  \end{table}

 Note that 4 and 5 in Table \ref{ti} were done in \cite{ds}.

\section{When  is a reduced module semisimple?}

It is known that  a semisimple module over a commutative ring is reduced. In this section, we investigate conditions 
 when the converse to this holds for modules over arbitrary rings. 

\begin{thm}\label{Th1} Let $R$ be a  ring and $M$  an $R$-module such that  either the ring $R$ or the ring $R/\text{ann}_R(M)$
is embeddable in the module $M$. If $M$ is artinian and reduced, then every $R/\text{ann}_R(M)$-module is    semisimple. 
\end{thm}

\begin{prf}
 If there exists a monomorphism from $R$ into $M$, then $R$ is isomorphic to a submodule of $M$. Since $M$ is artinian, so is $R$ and hence 
 $R/\text{ann}_R(M)$ is also artinian. On the other hand, if $R/\text{ann}_R(M)$ is embeddable in an artinian module $M$, 
 then $R/\text{ann}_R(M)$  is also artinian.
 It is easy to prove that if $M$  is a reduced $R$-module, then   $R/\text{ann}_R(M)$ is a reduced ring. 
 However, a reduced artinian ring is semisimple. 
 Hence,  $R/\text{ann}_R(M)$ is a semisimple ring. So,  every $R/\text{ann}_R(M)$-module is   semisimple.
\end{prf}

 \begin{prop}\label{p1}
  In the  $R$-modules $M$  given in 1-2, $R$ is embeddable in $M$; and in the $R$-modules  $M$ 
  given in 3-4, $R/\text{ann}_R(M)$ is embeddable in $M$.
  \begin{enumerate}
  \item Free $R$-modules $M$,
  \item $M=\oplus_iP$ where $P$ is a generator module for $R$-Mod,
  \item  finitely generated modules $M$ over commutative rings $R$,
  \item completely prime $R$-modules $M$.
 \end{enumerate}
 \end{prop}

\begin{prf}
 \begin{enumerate}
  \item It is enough to remember that if $M$ is a free $R$-module, then $M$ is isomorphic to the $R$-module $R^n$ with $n$ copies of $R$
  for some    positive integer $n$.
  \item A module   $P$  is a generator module  for $R$-Mod if and only if $R$ is a direct summand of a direct sum $M=\oplus_{i}P$,
  see \cite[Theorem 18.8]{Lam}. 
  \item Let $M$ be a finitely generated module over a commutative ring $R$ with generators $\{m_1, m_2, \cdots, m_n\}$.
  Define an $R$-homomorphism
  $f:R/\text{ann}_R(M) \rightarrow M$ by $f(\bar{r})= (\bar{r}m_1, \bar{r}m_2, \cdots, \bar{r}m_n)$.  $f$  is a monomorphism. If $f(\bar{r})=0$, 
  then $\bar{r}m_i=0$ for all
  $i\in\{1, 2, \cdots, n\}$ and  $\bar{r}M=\bar{r}\sum_{i=1}^{n}Rm_i= \sum_{i=1}^{n}R\bar{r}m_i=0.$ Thus, $\bar{r}\in \text{ann}_R(M)$. 
  \item  Let $f:R/\text{ann}_R(M)\rightarrow M$  be defined by $f(\bar{r})=rm$ where  $\bar{r}=r+\text{ann}_R(M)$.
  Suppose that $f(\bar{r})=0$, then $rm=0$. Since $M$ is completely prime,  
   $m=0$ or $rM=0$. If $m=0$, $\text{Ker}~f = R/\text{ann}_R(M)$ and $f$  is the zero homomorphism. If $rM=0$, $r\in \text{ann}_R(M)$
    and $\bar{r}=0$. So, $f$  is injective. In both cases, $R/\text{ann}_R(M)$ is embeddable in $M$.

 \end{enumerate}

\end{prf}
In Proposition \ref{p1}, we have given a desirable property that completely prime modules possess but prime  modules 
do not have; a yet another 
justification for completely prime modules. See \cite[Section 2]{S2} for other advantages completely prime modules  have over prime  modules.
In a prime $R$-module $M$, the ring $R/\text{ann}_R(M)$ need not be embeddable in $M$, see Example \ref{z}(2 \& 4) and Remark \ref{rma}(3).\\




\begin{cor} The following statements  are true: 
\begin{enumerate}
 \item A faithful artinian completely prime module is semisimple.
 \item An artinian, reduced and faithful free  module is semisimple.
 \item An artinian, reduced, faithful and finitely generated module over a commutative ring is semisimple.
 \item A reduced, artinian progenerator module is semisimple.
\end{enumerate}
 
\end{cor}

\begin{prf}
 It is enough for one to see that a completely prime module is reduced,  
 a progenerator module is a generator module,  and by  \cite[Remark 18.9(B)]{Lam} any generator module is faithful.
 The rest follows from Theorem \ref{Th1} and Proposition \ref{p1}.

\end{prf}

 Let $A_n$ be the $n$th Weyl algebra. A finitely generated $A_n$-module is called {\it holonomic} if it is zero, or it has dimension $n$.
 
 \begin{cor}\label{above} A faithful holonomic  module    is   reduced    if  and only if it is   semisimple and rigid.
  \end{cor}
 
 \begin{prf} 
  By \cite[Theorem 2.2]{C}, a holonomic module $_RM$ is artinian. By definition, a  holonomic module is a free module and hence 
  by Proposition \ref{p1},  $R$ is   embeddable in $M$. Applying Theorem \ref{Th1}  shows that $M$ is a semisimple $R$-module.
  The remaining part is easy since a reduced module  is always rigid. Conversely, if $M$ 
  is semisimple and rigid as an $R$-module, then by \cite[Corollary 2.30]{B} it  is reduced.
  \end{prf}
  
  \begin{cor} An artinian $R$-algebra $M$  is a reduced    $R/\text{ann}_R(M)$-module 
 if  and only if $M$ is a semisimple rigid   $R/\text{ann}_R(M)$-module.
  \end{cor}
  \begin{prf}
  An artinian $R$-algebra is an   artinian finitely generated module defined over a commutative  ring $R$. By Proposition \ref{p1}, $R/(0:M)$ is
  embeddable in $M$. The rest of the proof is {\it mutatis mutandis} the one given in Corollary \ref{above} above.
  \end{prf}


 


 \begin{cor}
  Let $(R, \mathfrak{m})$ be a Gorenstein artinian local ring. If the  injective envelope $\mathcal{E}_R(R/\mathfrak{m})$ of $R/\mathfrak{m}$ is 
  a reduced $R$-module, then $\mathcal{E}_R(R/\mathfrak{m})$ is   a semisimple $R$-module.
 \end{cor}

 \begin{prf}
  By \cite[Example 3.2]{Kikumichi}, $\mathcal{E}_R(R/\mathfrak{m})$  is a finitely generated faithful $R$-module in which $R$ embeds. Since 
  artinian modules are closed under taking injective envelopes, $\mathcal{E}_R(R/\mathfrak{m})$  is an artinian module.
  Thus, $\mathcal{E}_R(R/\mathfrak{m})$  is an artinian reduced faithful module in which $R$ embeds. By Theorem \ref{Th1},
  $\mathcal{E}_R(R/\mathfrak{m})$    is a semisimple $R$-module.
  
 \end{prf}



 A zero submodule of an $R$-module $M$  is said to {\it satisfy the complete radical formula} 
 
 (resp. {\it satisfy  the radical formula}) if  $\langle E_M(0)\rangle = \beta_{co}(M)$ (resp. $\langle E_M(0)\rangle = \beta(M)$) where 
 $\beta_{co}(M)$ (resp. $\beta(M)$) is the completely prime radical of $M$, i.e., the intersection of all completely prime submodules  of $M$
 (resp. the prime radical of $M$, i.e., the intersection of all prime submodules of $M$).
\begin{thm}\label{zero}
If $R$ is a left artinian  ring, and  the zero submodule of an $R$-module $M$  satisfies the complete radical formula; 
then $M$ is semisimple whenever it is reduced. 
\end{thm}

\begin{prf} Since $R$ is artinian and therefore $R/J(R)$ is also artinian, by \cite[Exercise  4. 16]{Lam2}, 
to show that $M$  is semisimple, it is enough to show that $J(R)M=0$  where $J(R)$ is the Jacobson
 radical of the ring $R$. 
 Suppose that $\langle E_M(0)\rangle = \beta_{co}(M)$.  If $M$  is reduced  $\mathcal{N}(M)=0$
  and $\langle E_M(0)\rangle =0$ by Proposition \ref{nil}. It follows that  $\beta_{co}(M)=0$.
 Since $\beta(R)M\subseteq \beta_{co}(R)M\subseteq \beta_{co}(M)$ (see \cite[Lemma 5.4]{DS}),  where $\beta_{co}(R)$ (resp. $\beta(R)$)
 is the completely prime radical (resp. prime radical) of $R$, we have 
   $\beta(R)M=0$. Since $R$ is left artinian, $\beta(R)=J(R)$ and hence $J(R)M=0$ as required.
\end{prf}

\begin{cor}\label{co}
Any module over  a commutative artinian  ring    is reduced if and only if it  is semisimple. 
\end{cor}

\begin{prf}
A semisimple module over a commutative ring  is reduced. 
Any commutative artinian ring $R$ satisfies the radical formula (see \cite{Man})
and hence the zero submodule of an $R$-module satisfies the radical
 formula. However,  for modules over commutative rings, there is no   distinction between modules that satisfy the radical
 formula and modules that satisfy the complete radical formula. The rest follows from Theorem \ref{zero}.
\end{prf}

\begin{cor}
A commutative artinian  ring  is reduced if and only if it  is semisimple. 
\end{cor}

\begin{prf}
 For any commutative ring $R$, $E_R(0)=\beta(R)$, i.e., the zero ideal of $R$ satisfies the radical formula. The rest follows from 
 Corollary  \ref{co}.
\end{prf}


  


 \section{Examples}
 
 \begin{enumerate}
  \item  If $R:=M_n(D)$ where $D$ is a division ring, then $_RR$ is semisimple  and  nil. This shows that in a module over a noncommutative 
   ring, existence  of nilpotent elements does not in general hinder semisimplicity. This is contrary to what happens for modules over 
   commutative rings.
  \item  If $\Z$ is the ring  of integers, then $_{\Z}\Z$ is reduced  but not semisimple.  
     
  \item \label{3} Let $Q: = $ \begin{tikzpicture}
   \draw[->, thick] (0, 2) -- (1.4, 2);
   \coordinate [label=left: .] (.) at (0.1,2);
\coordinate [label=right:. ] (.) at (1.3 ,2); 
   \end{tikzpicture}   
  be a quiver. Its path algebra $kQ$  has as basis $\{e_1, e_2, a\}$. Let $M$ be a 2-dimensional $kQ$-module with basis   $\{x, y\}$ and the  
  action of  $kQ$ on  $M$ be given by 
\begin{center}
 
\begin{tabular}{c|cc}
     & $x$ &  $ y$\cr\hline
  $e_1$ & $x$ & $0$\cr
  $e_2$ & $0$ & $y$ \cr
  $a$ & $0$ & $x$ \cr
\end{tabular}
\end{center}
$M$ has only one simple submodule  which is generated by $x$. So, $M$  is not semisimple. 
If $m_1=\mu x + \lambda y, m_2=\mu x$, $m_3=\lambda y\in M$ and 
$\alpha= pe_1+qe_2+ra\in kQ$, then $(ra)^2m_1=(ra)^2m_3=0$ but $(ra)m_1\not=0$ and $(ra)m_3\not=0$. This shows that $m_1$ and $m_3$ are
 nilpotent elements of $M$. However, $m_2$ is not nilpotent and hence $M$  is not reduced.  
 
 \item  If the action in Example \ref{3}  above is changed to
 \begin{center}
 
\begin{tabular}{c|cc}
     & $x$ &  $ y$\cr\hline
  $e_1$ & $x$ & $0$\cr
  $e_2$ & $0$ & $y$ \cr
  $a$ & $y$ & $x$ \cr
\end{tabular}
\end{center}
then we get a simple module and hence a semisimple module which is reduced.
 \end{enumerate}

 \begin{exam}[\rm{\bf The Golden Example}]  \label{z} \rm
 Let $M:=M_n(\Z/k\Z)$ where $k$  is an integer be the ring of all $n\times n$  matrices defined over the ring $\Z/k\Z$,
 $$N:=\left\{(m_{ij})\in M~:~ \sum_{j=1}^nm_{ij}=0~(\text{mod}~k) ~\forall i\in \{1, 2, \cdots, n\}\right\}~\text{and}~ 
 R:=M_n(\Z).$$ 
 \begin{enumerate}
  \item $M$ is an $R$-module with order $k^{n^2}$ and  $N$  is a minimal submodule of $M$  with order $k^n$;
  \item As an $R$-module, $N$  is simple. Hence, it is      prime   with $\text{ann}_R(N)=M_n(k\Z)$;
    \item  $N$  is not a completely prime module,  i.e., there exists a nonzero element $m\in N$ such that $\text{ann}_R(N)\subsetneq \text{ann}_R(m)$;
    \item The homomorphism $f: R/\text{ann}_R(N)\rightarrow N$ given by $f(\bar{r})=rn$ where $\bar{r}= r + \text{ann}_R(N)$ and $n\in N$ is not 
    injective;
    
    \item $N$  is a nil $R$-module;
    \item $0=\beta(M)\subsetneq \beta_{co}(M)=E_M(0)= \mathcal{N}(M)=M$.
    
 \end{enumerate}
 
 \end{exam}
 
 \begin{rem}\rm \label{rma}
 
 We call Example \ref{z} the  Golden Example because it is used for many purposes; we list them below:
 \begin{enumerate}
  \item It is an example of a nil module.
  \item It generalizes an example of a prime module   which is not completely prime  given  in \cite[Example 1.3]{S2} and of a classical 
         completely prime module which is not classical prime given in \cite[Example 3.1]{CCP}. In these two cases, $k$ and $n$ were
          taken to be $k=n=2$.
   \item  It shows that in a prime $R$-module $M$, the ring $R/\text{ann}_R(M)$ need not embed in $M$. More concretely, take $k=2$ in 
   Example \ref{z}. This gives another advantage completely prime modules have over  prime modules.
   \item It shows that a prime module can be nil. We  show in Section \ref{Sectionl} that $s$-prime, $l$-prime and completely prime modules
   cannot be nil.
   \item It is an example of a module that satisfies the complete radical formula but it is neither 2-primal nor satisfies the radical
    formula; see a particular case used in \cite{cp2p}.
    
    \item Since $\beta(M)\not=\beta_{co}(M)$ but $\beta_{co}(M)=\mathcal{N}(M)$, it gives a negative answer to  Question 8.2.1  posed 
     in \cite{thesis} as to whether we can succeed in writing 2-primal modules in terms of nilpotent elements as  is the case for 2-primal rings.
 \end{enumerate}
  \end{rem}

%

 \begin{exam}\rm
  Let $A$ be a ring which is not necessarily unital.
  Define $A^1:=\{(a, n)~:~a\in A, n\in \Z\}$ with component-wise addition and multiplication given by 
 $(a, n)(b, m)=(ab+am+nb, nm)$. Then $A^1$ is a ring (called the Dorroh extension of $A$) with unity $(0, 1)$. If $A$ is  such that 
 $(A, +)$  has no torsion elements and $A$ as a ring has at least one nonzero nilpotent element, then the submodule $B:=\{(0, n)~:~n\in \Z\}$ of
 the $A^1$-module $A^1$ is nil.
 For if $a$ is a nonzero nilpotent element of $A$, then $(a, 0)(0, n)=(an , 0)\not=(0, 0)$ for all $0\not=n\in \Z$ since $(A, +)$ 
 has no torsion elements. This shows that $(a, 0)A(0, n)\not=0$.
 However, since $a^k=0$ for some positive integer $k$, we have $(a, 0)^k(0, n)=(a^kn , 0)=(0, 0)$ which shows that $B$  is nil.

  \end{exam}

 \begin{exam}\rm

Let $R$  be  a commutative  ring  and $S:=R[x, y]$, the polynomial ring over $R$ in the indeterminates $x$ and $y$.
The regular module $_SS$ is reduced. However, a finite dimensional $S$-module
$R[x, y]/\langle x^4, xy^2, x^3y, y^4\rangle$ has some nilpotent elements.
The generators of $R[x, y]/\langle x^4, xy^2, x^3y, y^4\rangle$  are 
$\{1, x, x^2, x^3, y, y^2, y^3, xy, x^2y\}$  all nilpotent except $\{y^3, x^2y, x^3\}$.
This information is encoded in a combinatorial object given in Figure \ref{fig}. 
 
 \begin{figure}[h]
 \begin{center}
  \begin{tikzpicture}
   \draw[->, thick] (5, 2) -- (10.4, 2); 
   \draw[->, thick] (5, 2) -- (5, 5.4);
   \draw[-] (5, 2.6) -- (7.4, 2.6);
   \draw[-] (5, 3.2) -- (6.8, 3.2);
   \draw[-] (5, 3.8) -- (5.6, 3.8);
   \draw[-] (5, 4.4) -- (5.6, 4.4);
   \draw[-] (5.6, 2) -- (5.6, 4.4);
   \draw[-] (6.2, 2) -- (6.2, 3.2);
   \draw[-] (6.8, 2) -- (6.8, 3.2);
   \draw[-] (7.4, 2) -- (7.4, 2.6);
   \draw[red, thick] (5.3, 4.1)  circle (0.32cm);
    \draw[red, thick] (6.5, 2.9)  circle (0.32cm);
    \draw[red, thick] (7.1, 2.3)  circle (0.32cm);
   
   \coordinate [label=left: $x^3$] (.) at (7.5,2.2);
   \coordinate [label=left: $x^2$] (.) at (6.9,2.2);
   \coordinate [label=left: $x$] (.) at (6.2,2.2);
   \coordinate [label=left: $1$] (.) at (5.5,2.2);
   \coordinate [label=left: $y$] (.) at (5.5,2.8);
   \coordinate [label=left: $y^2$] (.) at (5.6,3.4);
   \coordinate [label=left: $xy$] (.) at (6.3,2.8);
   \coordinate [label=left: $x^2y$] (.) at (7,2.8);
   \coordinate [label=left: $y^3$] (.) at (5.6, 4);   
   \end{tikzpicture}   
  \caption{Generators of the $S$-module $R[x, y]/\langle x^4, xy^2, x^3y, y^4\rangle$.}
  \label{fig}
 \end{center}
 
 \end{figure}
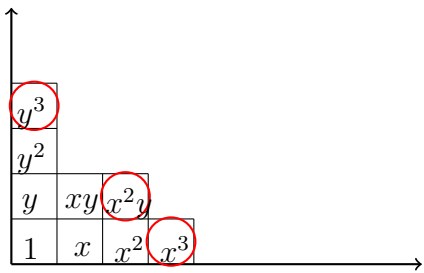
 Note that the non-nilpotent module elements
$\{y^3, x^2y, x^3\}$ have been circled and occur at the ``sharp points of the stairs'' in Figure \ref{fig}. In general, they can 
be determined easily by just drawing such combinatorial objects. We also remark that the same circled elements have connections to
Hilbert schemes. These connections are   outlined in  \cite[Chapter 7]{Hi}.
 
\end{exam}

\begin{exam}\rm 
 For any positive integer $n$ and any prime  integer $p$, $M_n(\Z/p\Z)$ is a reduced $\Z$-module and the $\Z$-module 
 $M_n(\Z/p^k\Z)$ for any positive integer $k$ greater than $1$ contains nilpotent elements.
\end{exam}

\section{Effect of nilpotents on primeness}\label{Sectionl}

 Andrunakievich in \cite{Andrunakievich} defined  $l$-prime and $s$-prime modules which were further studied in \cite{Levitzki} and 
 \cite{Kothe} respectively. We give the effect   nilpotent elements have 
 on $l$-prime and $s$-prime modules.  Let $\mathcal{L}(R)$ and $\mathcal{U}(R)$
 denote the Levitzki radical and upper nil radical respectively of a ring $R$.
 A module $M$  is said to be {\it $l$-prime} (resp. {\it $s$-prime}) if it is prime and  $\mathcal{L}(R/\text{ann}_R(M))=0$
 (resp. if it is prime and  $\mathcal{U}(R/\text{ann}_R(M))=0$), see \cite{Levitzki} (resp. \cite{Kothe}) for the other equivalent 
 definitions.

  \begin{prop}\label{7}
  Let $M$  be an $R$-module. If any one of the following is true, then $\mathcal{N}(M)\not=0$, i.e., $M$  contains a nilpotent element.
  \begin{enumerate}
   \item $\mathcal{U}(R/\text{ann}_R(M))\not=0$,
    \item $\mathcal{L}(R/\text{ann}_R(M))\not=0$,
     \item $\mathcal{\beta}(R/\text{ann}_R(M))\not=0$.
  \end{enumerate}
  \end{prop}

  \begin{prf}
   Suppose  $\mathcal{U}(R/\text{ann}_R(M))\not=0$. Then there exists $0\not= I/\text{ann}_R(M)\vartriangleleft R/\text{ann}_R(M)$
    such that $I/\text{ann}_R(M)$ is nil. Then for all $\bar{0}\not=\bar{r}\in I/\text{ann}_R(M)$, there exists a positive integer $k$ such 
    that $\bar{r}^k=\bar{0}$. This is equivalent to saying that, for all $r\in I\setminus \text{ann}_R(M)$, there exists a positive 
    integer $k$ such that $r^k\in \text{ann}_R(M)$. So, $rM\not=0$ and $r^kM=0$. This implies that there exists a nonzero element $m$  in 
    $M$ such that $r^km=0$ and $rm\not=0$ so that $rRm\not=0$. This shows that  $m\in \mathcal{N}(M)$. Since 
    $\mathcal{\beta}\subseteq \mathcal{L}\subseteq \mathcal{U}$, if either 
    $\mathcal{L}(R/\text{ann}_R(M))\not=0$ or $\mathcal{\beta}(R/\text{ann}_R(M))\not=0$, it follows that, $\mathcal{U}(R/\text{ann}_R(M))\not=0$,
    a case which is already proved.
  \end{prf}

  In Corollary \ref{f} below, we retrieve an already known result.
  
  \begin{cor}\label{f}
   If $M$  is a prime and reduced $R$-module, then 
   \begin{enumerate}
    \item $M$  is $s$-prime,
    \item $M$  is $l$-prime.
       \end{enumerate}

  \end{cor}
\begin{prf}
 If $M$  is reduced, $\mathcal{N}(M)=0$. By Proposition \ref{7}, $$\mathcal{U}(R/\text{ann}_R(M))=\mathcal{L}(R/\text{ann}_R(M))=0.$$
 Since $M$  is in addition prime, it follows from \cite[Corollary 2.1]{Kothe} (resp. \cite[Proposition 2.2]{Levitzki})  that 
 $M$  is $s$-prime (resp. $l$-prime).
\end{prf}
 This result is already known, because from \cite[Theorem 2.10 $\&$ Remark 2.11]{DS} and \cite[Corollary 2.4]{2p}
 a prime and reduced module is completely prime. 
 It was shown in \cite[Propositions 3.1 $\&$ 3.2]{2p}  and \cite[Proposition 2.7]{Levitzki} 
 that  a completely prime module is $s$-prime and an $s$-prime module is $l$-prime   respectively.\\
 
 Corollary \ref{VC} shows that ``reduced'' in a faithful $R$-module carries over  in some sense to the ring $R$.
 
 \begin{cor}\label{VC}
  If $M$ is a faithful and reduced $R$-module, then 
  \begin{enumerate}
   \item $\mathcal{U}(R)=0$, i.e., $R$ has no nonzero nil ideals;
   \item $\mathcal{L}(R)=0$, i.e., $R$ has no nonzero locally nilpotent ideals;
   \item $\mathcal{\beta}(R)=0$, i.e., $R$ has no nonzero nilpotent ideals.
  \end{enumerate}

 \end{cor}
 
 \begin{cor}
  A prime module which is not $l$-prime, not $s$-prime or not completely prime, contains nonzero nilpotent elements.
 \end{cor}

  It is impossible to write a nil module as a sum of completely prime modules. Hence, it is impossible
   to have a nil module defined over  a commutative ring written  as a sum of prime  modules. If it were possible, the completely prime 
   modules making up the sum would be nil which is a contradiction - completely prime modules are always reduced. It is important also
   to note that in a module over a commutative ring, prime modules are completely prime. So, existence of nilpotents inhibits some structure;
    the structure of having a sum of completely prime modules.

\begin{cor}
 For any module $M$, $\mathcal{N}(M)=0$  implies that  $$\mathcal{U}(M)=\mathcal{L}(M)=\beta(M)=0.$$
\end{cor}

\begin{prf}
 If $\mathcal{N}(M)=0$, then $M$  is reduced and $\beta_{co}(M)=0$. This leads to the desired result  since 
 $\beta(M)\subseteq \mathcal{L}(M)\subseteq \mathcal{U}(M)\subseteq \beta_{co}(M)$.
\end{prf}

For rings, the sum of all nil ideals of a ring coincides with the intersection of all $s$-prime ideals of that ring. We get the following
 question.
 
 \begin{qn}{\rm \cite[Question 8.2.3]{thesis}}\label{q}
  How does the upper nil radical of a module $M$ which is given as the  intersection of all $s$-prime submodules of $M$ compare 
  with the sum of all nil submodules of $M$?
 \end{qn}

 In conclusion, nilpotent elements of a module control its structure. They inhibit semisimplicity for modules defined over commutative rings.
They do not allow a module to be any of the following: torsion-free, completely prime, $l$-prime and $s$-prime. They determine
in addition to some other conditions whether  a module can be written as a direct sum of prime submodules or not. In a situation where they
do not appear, i.e., when the module is reduced,  the module behaves nicely; for instance, it behaves as though it is  defined 
over a commutative ring, i.e., it has the insertion of factor property, it is 2-primal and it is symmetric, see \cite[Theorems 2.2 \& 2.3]{2p}.
  \section*{Acknowledgement}
  
  I wish to thank professors Ken A. Brown and Michael Wemyss both of University of Glasgow, and Professor Balazs Szendroi of 
  University of Oxford for  the hospitality during my visit
  to the UK and for the discussions that greatly improved this work. This research was supported by LMS and 
  Sida bilateral programme (2015-2020) with Makerere University; project 316: Capacity building in Mathematics and its Applications.

\end{document}